\documentclass[12pt]{article}
\usepackage{amsfonts}
\usepackage{amsmath}
\usepackage{amssymb,amsmath}
\usepackage{amscd}
\newfont{\bbb}{cmb10 at 10pt}

\def\a{\alpha}
\def\G{\Gamma}

\def\b{\beta}
\def\e{\varepsilon}

\def\cl{\centerline}

\def\vs{\vspace*}

\def\Z{\mathbb{Z}}

\def\F{\mathbb{F}}
\def\QED{\hfill$\Box$}

\numberwithin{equation}{section}
\newtheorem{theo}{Theorem}[section]
\newtheorem{defi}[theo]{Definition}
\newtheorem{coro}[theo]{Corollary}
\newtheorem{lemm}[theo]{Lemma}
\newtheorem{prop}[theo]{Proposition}

\newtheorem{rema}[theo]{Remark}
\newtheorem{ex}[theo]{Example}

\begin{document}

\cl{{\bf Hom-Lie color algebra structures }}\vs{6pt}

\cl{Lamei Yuan}

\cl{\small Department of Mathematics, \!University of Science  and
Technology \!of \!China, Hefei 230026, P. R. China}

\cl{\small E-mail: lmyuan@mail.ustc.edu.cn,
yuanlamei123@126.com}\vs{6pt}

{\small\parskip .005 truein \baselineskip 3pt \lineskip 3pt

\normalsize\noindent{\bbb Abstract} The aim of this paper is to
introduce the notion of Hom-Lie color algebras. This class of
algebras is a natural generalization of the Hom-Lie algebras as well
as a special case of the quasi-hom-Lie algebras. In the paper,
homomorphism relation between Hom-Lie color algebras are defined and
studied. We present a way to obtain Hom-Lie color algebras from the
classical Lie color algebras along with algebra endomorphisms and
offer some applications. Besides, we introduce a multiplier $\sigma$
on the abelian group $\G$ and provide constructions of new Hom-Lie
color algebras from old ones by the $\sigma$-twists. Finally, we
explore some general classes of Hom-Lie color admissible algebras
and describe all these classes via $G$-Hom-associative color
algebras, where $G$ is a subgroup of the symmetric group $S_3$.

\vs{5pt}

\noindent{\bbb Key words:} Hom-Lie algebras, Hom-Lie color algebras,
Hom-Lie color admissible algebras, $G$-Hom-associative color
algebras, Homomorphism, $\sigma$-twists}\\ {\bbb 2000 Mathematics
Subject Classification }\quad 17B75,
17B99\\

\parskip .001 truein\baselineskip 6pt \lineskip 6pt

\vs{18pt}

\cl{\bf\S1. \
Introduction}\setcounter{section}{1}\setcounter{equation}{0}

\vs{8pt} The motivations to study Hom-Lie structures are related to
physics and deformations of Lie algebras, in particular Lie algebras
of vector fields. The paradigmatic examples are $q$-deformations of
Witt and Virasoro algebras constructed in pioneering works (e.g.,
\cite{AS1, M, MPJ, MZP, JDS, K}). This kind of algebraic structures
were initially introduced by Hartwig, Larsson and Silvestrov
\cite{JDS} during the process of investigating the $q$-deformation
of Lie algebras. Later, it was further extended by Larsson and
Silvestrov to quasi-hom-Lie algebras and quasi-Lie algebras
\cite{DS1,DS2}, while the Hom-Lie algebra structures were more
detailed studied in \cite{AS}, including the Hom-Lie admissible
algebras and $G-$Hom-associative algebras, both of which can be
viewed as generalizations of Lie admissible algebras and
$G-$associative algebras, respectively. Recently, this kind of
algebras was considered in $\Z_2$-graded case by Ammar and Makhlouf
\cite{FA} and thus is said to be Hom-Lie superalgebras. The main
feature of quasi-Lie algebras, quasi-hom-Lie algebras and Hom-Lie
(super)algebras is that the skew-symmetry and the Jacobi identity
are twisted by several deforming twisting maps, which lead to many
interesting results.\par

 The Lie admissible algebras were introduced by A. A. Albert in 1948 \cite{AA}.
Physicists attempted to introduce this structure instead of Lie
algebra. For instance, the validity of Lie-admissible algebras for
free particles is well known. These algebras arise also in classical
quantum mechanics as a generalization of conventional mechanics (see
\cite{HC2,HC3}). The study of flexible Lie admissible algebra was
also initiated in [\ref{AA}] and has investigated in a number of
papers (see e.g., \cite{ME, HC}). The authors in [\ref{AS}] extended
to Hom-Lie algebra the classical concept of Lie admissible algebras,
while Hom-Lie admissible superalgebras were considered in \cite{FA}.

In the present paper, we introduce and study the Hom-Lie color
algebras, which can be viewed as an extension of Hom-Lie
(super)algebras to $\G$-graded algebras, where $\G$ is any abelian
group.
 This kind of algebras is a special case of the quasi-hom-Lie algebras from \cite{DS1} and contains as a subclass on the one hand the
 Lie color algebras and in particular Lie (super)algebras, and on
 the other hand various known and new single and multi-parameter
 families of algebras obtained using algebra endomorphisms and twisted
 mappings of
 Lie color algebras and algebras
 of vector fields.

 The paper is organized as follows. In Sec.$2$  we collect our
 conventions and consider some generalities on graded algebra, Lie color
 algebra and (quasi-)Hom-Lie algebra structures.

Hom-Lie color algebras are defined in Sec.$3$ and examples of such
kind of algebras are given. Homomorphism between two Hom-Lie color
algebras are also defined and studied. In particular, we show that
how arbitrary Lie color algebras deform into Hom-Lie color algebras
via algebra endomorphisms and offer some applications.

In Sec.$4$ we introduce a multiplier $\sigma$ on the abelian group
$\G$ and provide constructions of new Hom-Lie color algebras using
the twisting action of the multiplier $\sigma$. We show that the
$\sigma$-twist of any Hom-Lie color algebra is still a Hom-Lie color
algebra.

In Sec.$5$ we extend the classical concept of Lie admissible
algebras to Hom-Lie color settings. Hence, we obtain a more
generalized algebra class called Hom-Lie color admissible algebras.
We also explore some other general class of algebras:
$G-$Hom-associative color algebras, where $G$ is any subgroup of the
symmetric group $S_3$, using which we classify all the Hom-Lie color
admissible algebras.

Throughout this work, $\mathbb{F}$ denotes a field of characteristic
zero and $\Gamma$ stands for an abelian group. The multiplicative
group of nonzero element of $\mathbb{F}$ is denoted by
$\mathbb{F}^*$. All vector spaces and algebras are assumed to be
over $\mathbb{F}$. All gradations are understood to be
$\G$-gradations.

\vs{18pt}

\cl{\bf\S2. \
Preliminaries}\setcounter{section}{2}\setcounter{equation}{0}
\vs{8pt}

 Let us begin with
some definitions concerning graded algebraic structures. For a
detailed discussion of this subject we refer the reader to the
literatures (e.g. \cite{M1} and references therein).

 Let $\Gamma$ be an abelian group.  A vector space $V$ is said to be $\Gamma$-graded, if there is a family $(V_{\alpha})_{\alpha\in\Gamma}$ of subspaces of $V$ such that
$V=\bigoplus _{\alpha\in \Gamma} V_{\alpha}$.
 An element $a\in V$ is said to be
homogeneous of degree $\alpha$ if $a\in V_{\alpha}$, $\alpha\in
\Gamma$, and in this case, $\alpha$ is called the color of $a.$ As
usual, denote by \={a} the color of an element $a\in V$. Thus each
homogeneous element $a$ in $V$ determines a unique group element
$\bar a\in\G$ by $a\in V_{\bar a}$. Fortunately, we can almost drop
the symbol $``-"$, since confusion rarely occurs. We will denote by
$\mathcal {H}(V)$ the set of all the homogeneous elements of $V$ in
the sequel.
\par

Let $V=\oplus_{\alpha\in\Gamma}V_{\alpha}$ and
$W=\oplus_{\alpha\in\Gamma}W_{\alpha}$ be two $\Gamma$-graded vector
spaces. A linear mapping $f:V\rightarrow W$ is said to be
homogeneous of degree $\gamma\in\Gamma$ if
$$f(V_{\alpha})\subset W_{\alpha+\gamma},\ \ \ \ \forall\  \  \alpha\in\Gamma.$$
If, in addition, $f$ is homogenous of degree zero, namely,
$f(V_{\alpha})\subset W_{\alpha}$ holds for any $\alpha\in\Gamma$,
then $f$ is said to be even.\par

An algebra $A$ is said to be $\Gamma$-graded if its underlying
vector space is $\Gamma$-graded, i.e.
$A=\oplus_{\alpha\in\Gamma}A_{\alpha}$, and if, furthermore,
$A_{\alpha}A_{\beta}\subset A_{\alpha+\beta}$, for all $\alpha,\
\beta\in\Gamma.$ It is easy to see that if $A$ has a unit element
$e$, it follows $e\in A_0.$ A subalgebra of $A$ is said to be graded
if it is graded as a subspace of $A$. \par

Let $B$ be a second $\Gamma$-graded algebra. A homomorphism
$f:A\rightarrow B$ of $\Gamma$-graded algebras is by definition a
homomorphism of the algebra $A$ into the algebra $B$, which is, in
addition, an even mapping. \vs{6pt}

 The following definition is well-known from the theory of graded
 algebra.
\begin{defi}\label{def1}
Let $\G$ be an abelian group. A bi-character on $\G$ is a map
$\varepsilon :\ \G \times \G\rightarrow \F^*$ satisfying
$$
\varepsilon(\alpha,\beta)\varepsilon(\beta,\alpha )=1,$$
$$\varepsilon(\alpha,\beta+ \gamma)=\varepsilon(\alpha ,\beta)\varepsilon( \alpha,\gamma
 ),$$
$$\varepsilon(\alpha+\beta,\gamma)=\varepsilon(\alpha ,\gamma)\varepsilon( \beta,\gamma
 ),$$
for all $\alpha, \beta,\gamma \in \G.$
\end{defi}

It is easy to see  that
$\varepsilon(\alpha,0)=\varepsilon(0,\alpha)=1$ and $
\varepsilon(\alpha,\alpha)=\pm1,$ $\forall$ $\alpha\in \G.$ In
particular, for fixed  $\alpha\in\Gamma$, the induced mapping
$\varepsilon_{\alpha}:\Gamma\rightarrow \mathbb{F}^*$ defines a
homomorphism of groups by
$\varepsilon_{\alpha}(\beta)=\varepsilon(\alpha, \beta).$

\begin{defi}\label{def2} A $\G$-graded algebra $L=\oplus_{\a\in\G}L_\a$,
whose even bilinear product mapping is denoted by a pointed bracket
$[\cdot,\cdot]$, is called a ($\G$-graded) $\e$ Lie algebra (or Lie
color algebra) if the
following identities are satisfied: \\

$\bullet$ $[x,y]=-\e(x,y)[y,x]$ \ \ ($\e$-skew symmetry),\\

$\bullet$ $\circlearrowleft_{x,y,z}\e(z,x)[x,[y,z]\,]=0$ \ \
($\e$-Jacobi identity),\\for all the homogenous elements $x,y,z\in
L,$ and where $\circlearrowleft_{x,y,z}$ denotes summation over the
cyclic permutation on $x,y,z.$
\end{defi}

Note that we have abbreviated $\e(\bar{x},\bar{y})$ by $\e(x,y)$ in
the definition above.

\begin{defi}\label{def3} Suppose that $L=\bigoplus_{\a\in \G}L_{\a}$ and $L'=\bigoplus_{\a\in \G}L_{\a}'$
are two Lie color algebras. A linear map $f$ from $L$ to $L'$ is
called an even linear map, if $f(L_\a)\subseteq L_\a'$ holds for any
$\a\in \G$. If, in addition, it satisfies
$f([x,y])=[f(x),f(y)]^{'}$, for all $x, \ y\in L$, then it is said
to be a Lie color algebra homomorphism.
\end{defi}

The Hom-Lie algebras were initially introduced by Hartwig, Larsson
and Silvestrov in \cite{JDS} motivated by examples of deformed Lie
algebras coming from twisted discretizations of vector fields. This
kind of algebras includes Lie algebras as a subclass. Quasi-hom-Lie
algebras and Hom-Lie superalgebras  were introduced in \cite{DS2}
and \cite{FA}, respectively. Hom-Lie superalgebras are a special
case of quasi-hom-Lie algebras. Both of these two classes of
algebras can be viewed as extensions of Hom-Lie algebras.
Hom-associative algebras were introduced in \cite{AS}. Let us now
recall these definitions.

\begin{defi}
A {\rm Hom}-associative algebra is a triple $(V, \mu, \zeta)$
consisting of a linear space $V$, a bilinear map $¦Ì\mu : V \times V
\rightarrow V $and an algebra homomorphism $\zeta: V \rightarrow V$
satisfying $$ \mu\big(\zeta(x), \mu(y, z)\big) = \mu\big(\mu(x, y),
\zeta(z)\big).$$
\end{defi}

\begin{defi} A Hom-Lie algebra is a triple $(L,[\cdot,\cdot],\zeta)$
consisting of a vector space $L$, a bilinear mapping $[
\cdot,\cdot]:L\rightarrow L$, and an algebra homomorphism
$\zeta:L\rightarrow L$ satisfying

$\bullet$ $[x,y]=-[y,x]$ (skew-symmetry),

$\bullet$ $\circlearrowleft_{x,y,z}[\zeta(x),[y,z]\,]=0$ (Hom-Jacobi
identity),\\
for all $x,y,z\in L.$
\end{defi}
\begin{rema}\rm In some papers, such as \cite{AS}, $\zeta$ is only
required to be a linear map.
\end{rema}

\begin{defi}\label{d1} A quasi-hom-Lie algebra (or shortened by  qhl-algebra) is a tuple
$(L,[\cdot,\cdot],\alpha,\beta,\omega)$ where

$\bullet$ $L$ is a vector space,

$\bullet$ $[\cdot,\cdot]:\ L\rightarrow L$ is a bilinear map called
a bracket or product in $L,$

$\bullet$ $\alpha,\beta:\ L\rightarrow L$ are linear maps,

$\bullet$ $\omega: D_{\omega}\rightarrow \mathcal {L}_{\F}(L)$ is a
map with domain of definition $D_\omega\subseteq L\times L,$ and
where $\mathcal {L}_{\F}(L)$ denotes the linear space of
 $\F$-linear maps of $L$,\\ such that the following conditions
hold:\\

$\bullet$ ($\omega$-symmetry) The bracket satisfies a generalized
skew-symmetry
$$[x,y]=\omega(x,y)[y,x],\ \ \mbox{for\ all}\ (x,y)\in D_\omega;$$

$\bullet$ ($\beta$-twisting) The map $\a$ is a $\beta$-twisted
algebra homomorphism, i.e.,
$$[\a(x),\a(y)]=\b\circ \a[x,y],\ \ \mbox{for\ all}\ (x,y)\in L;$$

$\bullet$ (qhl-Jacobi identity) The bracket satisfies a generalized
Jacobi identity
$$\circlearrowleft_{x,y,z}\Big\{\omega(z,x)\big([\a(x),[y,z]\,]+\b[x,[y,z]\,]\big)\Big\}=0,$$

for all $(z,x), (x,y), (y,z)\in D_\omega.$
\end{defi}

By taking $\b={\rm id}_L$  and $\omega=-{\rm id}_{L\times L}$ in
Definition \ref{d1}, one gets Hom-Lie algebras. When $L$ is a
$\Z_2$-graded vector space, the class of Hom-Lie superalgebras
defined in \cite{FA} is obtained by setting $\b={\rm id}_L$ and
$\omega(x,y)=-(-1)^{\bar x\bar y}$, for all  $x,y\in \mathcal
{H}(L).$

\vs{8pt} \cl{\bf\S3. \ Hom-Lie color
algebra}\setcounter{section}{3}\setcounter{equation}{0} \vs{8pt}

In this section,  we extend Hom-Lie algebras to any $\Gamma$-graded
vector space, which leads us to Hom-Lie color algebras. We present a
way to obtain such kind of algebras from Lie color algebras and also
offer some applications of this method. Before the definition, let
us firstly introduce Hom-color algebras.

\begin{defi}\label{def} A {\rm Hom}-color algebra is a triple $(A, \mu, \zeta) $ consisting of a $\Gamma$-graded linear space $A$, an even bilinear map $¦Ì\mu : A\times A \rightarrow A $ and an
even homomorphism $\zeta: A \rightarrow A$ satisfying $$
\mu\big(\zeta(x), \mu(y, z)\big) = \mu\big(\mu(x, y),
\zeta(z)\big).$$
\end{defi}

By a homomorphism of Hom-color algebras $ f :
(A,\mu,\zeta)\rightarrow (A',\mu',\zeta')$ we mean an algebra
homomorphism from $A$ to $A'$ such that $f\circ\zeta =\zeta'\circ
f$, or, in other words, such that the following diagram

\[\begin{CD}
A@>f>>A'\\
@V\zeta VV      @VV\zeta'V\\
A@>>f>A'\\
\end{CD}\]
commutes.\par
 The following result says that any
$\Gamma$-graded  associative algebra can deform into a Hom-color
algebra along with any even linear self-map.
\begin{theo}\label{th1} Let $(A,\mu)$ be a $\Gamma$-graded
associative algebra and $\zeta: A\rightarrow A$ be an even linear
map such that
$$\zeta\circ \mu=\mu\circ\zeta^{\otimes 2}.$$
Then $(A, \mu_{\zeta}=\zeta\circ\mu, \zeta)$ is a Hom-color algebra.
Moreover, one has
$\zeta\circ\mu_{\zeta}=\mu_{\zeta}\circ\zeta^{\otimes 2}.$
\par Suppose that $(B,\mu{'})$ is another $\Gamma$-graded associative algebra
and that $\zeta{'}:B\rightarrow B$ is an even linear map such that
$\zeta{'}\circ\mu{'}=\mu{'}\circ\zeta'^{\otimes 2}$. If
$f:A\rightarrow B$ is an associative algebra homomorphism such that
$f\circ \zeta=\zeta'\circ f$, then
$f:(A,\mu_{\zeta},\zeta)\rightarrow
(B,\mu{'}_{\zeta{'}}=\zeta{'}\circ\mu{'},\zeta{'})$ is also a
homomorphism of Hom-color algebras.
\end{theo}
\noindent{\it Proof}.\ \ Using the hypothesis that $\zeta$ is an
even linear map, one has for any $x,y,z$ in $\mathcal {H}(A)$
\begin{eqnarray*}
\mu_{\zeta}\big(\mu_{\zeta}(x,y),\zeta(z)\big)&=&\zeta\circ\mu\big(\zeta\circ\mu(x,y),\zeta(z)\big)
=\zeta^2\mu\big(\mu(x,y),z\big)\\&=&\zeta^2\mu\big(x,\mu(y,z)\big)=\zeta\circ\mu\big(\zeta(x),\zeta\circ\mu(y,z)\big)\\
&=&\mu_{\zeta}\big(\zeta(x),\mu_{\zeta}(y,z)\big),
\end{eqnarray*}
 from which it follows
that  $(A, \mu_{\zeta}=\zeta\circ\mu, \zeta)$ is a Hom-color
algebra.
\par Furthermore, observe that both $\zeta\circ\mu_{\zeta}$ and
$\mu_{\zeta}\circ\zeta^{\otimes 2}$ are equal to
$\zeta\circ\mu\circ\zeta^{\otimes 2}.$\par Finally,  since
$f:A\rightarrow B$ is an algebra homomorphism that satisfies $f\circ
\zeta=\zeta'\circ f$ by hypothesis, then for all $x,y\in \mathcal
{H}(A)$ we have
$$f\circ\mu_{\zeta}(x,y)=f\circ\zeta\circ\mu(x,y)=\zeta'\circ f\circ\mu(x,y)=\zeta'\circ \mu'(f(x),f(y))=\mu'_{\zeta'}\circ f^{\otimes 2}(x,y).$$
This proves $f$ is a homomorphism of Hom-color algebras.\QED
\vspace{2mm}\par In the definition of quasi-hom-Lie algebra, taking
$L$ to be a $\G$-graded vector space, $\beta={\rm id}_L$, and
$\omega(x,y)=-\varepsilon(\bar{x},\bar{y})$ for some bi-character
$\varepsilon$ on $\G$, we get the following class of algebras.

\begin{defi}\label{def8}A Hom-Lie color algebra is a quadruple $\big(L,[\cdot,\cdot], \zeta,
\varepsilon\big)$ consisting of a $\Gamma$-graded space $L$, an even
bilinear mapping $[\cdot,\cdot]:L\times L\rightarrow L$, a
homomorphism $\zeta$ and a bi-character $\e$ on $\Gamma$ satisfying

$\bullet$ $[x,y]=-\e(x,y)[y,x]$ ($\e$-skew-symmetry),

$\bullet$ $\circlearrowleft_{x,y,z}\e(z,x)[\zeta(x),[y,z]\,]=0$ (Hom
$\e$-Jacobi identity),\\
 for all homogenous elements $x,y,z$ in $L$.
\end{defi}
\begin{ex} \rm It is clear that Lie color algebras are examples of Hom-Lie color algebras by setting $\zeta={\rm id}_L$. If, in addition, $\e(x,y)=1$ or
$\e(x,y)=(-1)^{\bar{x}\bar{y}}$, then the Hom-Lie color algebra is
nothing but a classical Lie algebra or Lie superalgebra. The Hom-Lie
algebra and Hom-Lie superalgebra are also obtained when $\e(x,y)=1$
and $\e(x,y)=(-1)^{\bar{x}\bar{y}}$, respectively.
\end{ex}
\begin{ex}\label{ex1}\rm Let $\big(L,[\cdot,\cdot]_L,\zeta_L,\e\big)$ be a Hom-Lie
color algebra. Then the vector space $L':=L\otimes
\mathbb{F}[t,t^{-1}]$ can be considered as the algebra of Laurent
polynomials with coefficients in the Hom-Lie color algebra $L$. Note
that $L'$ can be endowed with a natural $\G$-grading as follows: an
element $x\in L'$ is said to be homogeneous of degree $\a\in\G$, if
there exist an element $x_{\a}\in L$ with degree $\a$
 and $f(t)\in\F[t,t^{-1}]$, such that $x=x_{\a}\otimes f(t).$
 Put $\zeta_{L'}=\zeta_L\otimes \rm id$ and define an even bilinear multiplication $[\cdot,\cdot]_{L'}$ on $L'$ by
$$[x\otimes f(t),y\otimes g(t)]_{L'}=[x,y]_L\otimes f(t)g(t),$$
for all $x,y\in \mathcal {H}(L)$ and $f(t),g(t)\in \F[t,t^{-1}]$.
With these definitions
$\big(L',[\cdot,\cdot]_{L'},\zeta_{L'},\e\big)$ is a Hom-Lie color
algebra. The verification of this consists of checking the axioms
from Definition \ref{def8} of Hom-Lie color algebras. For any
homogeneous elements $x,y,z\in L$, and $f(t),g(t),h(t)\in
\F[t,t^{-1}]$, it follows:
$$[x\otimes f(t),y\otimes g(t)]_{L'}=[x,y]_L\otimes f(t)g(t)=-\e(x,y)[y,x]_L\otimes g(t)f(t)=-\e(x,y)[y\otimes g(t),x\otimes f(t)]_{L'},$$
and
\begin{eqnarray*}
\circlearrowleft_{x,y,z}\e(z,x)[\zeta_{L'}(x\otimes f(t)),[y\otimes
g(t),z\otimes
h(t)]_{L'}]_{L'}&=&\circlearrowleft_{x,y,z}\e(z,x)[\zeta_L(x)\otimes
f(t),[y,z]_L\otimes g(t)h(t)]_{L'} \\
&=&\circlearrowleft_{x,y,z}\e(z,x)[\zeta_L(x),[y,z]_{L}]_{L}\otimes
f(t)g(t)h(t)\\&=&0,
\end{eqnarray*}
since $(L,[\cdot,\cdot]_{L},\zeta_L,\e)$ is a Hom-Lie color algebra.
\end{ex}

It is well known that there is always a Lie algebra associated to an
associative algebra via the commutator bracket. A Hom-associative
algebra can also gives rise to a Hom-Lie (super)algebra via the
commutator bracket \cite{AS,FA}. The following proposition presents
a similar result. In this sense, Hom-color algebras play the role of
associative algebras in the Hom-Lie color setting.

\begin{prop}\label{p1} Let $A=\bigoplus_{\alpha\in \Gamma}A_{\alpha}$ be a $\Gamma$-graded
vector space, $\e$ be a bi-character on $\Gamma$ and $(A, \mu,
\zeta)$ be a Hom-color algebra. One can define the color-commutator
on homogeneous elements $x,y\in A$ by $$[x,y] = \mu(x,y)
-\varepsilon(x,y)\mu(y,x),$$ and then extending by linearity to all
elements in $A$. Then $\big(A, [\cdot,\cdot],\zeta,\e\big)$ is a
Hom-Lie color algebra.
\end{prop}
\noindent{\it Proof}.\ \  The bracket is obviously
$\e$-skew-symmetric and the Hom-$\e$-Jacobi identity can be obtained
by a direct computation. For any homogeneous elements $x,y,z\in A,$
we have
\begin{eqnarray*}
\e(z,x)[\zeta(x),[y,z]\,]&=&\e(z,x)\mu\big(\zeta(x),\mu(y,z)\big)-\e(x,y)\mu\big(\mu(y,z),\zeta(x)\big)\\&-&\e(z,x)\e(y,z)\mu\big(\zeta(x),\mu(z,y)\big)+\e(y,z)\e(x,y)\mu\big(\mu(z,y),\zeta(x)\big),\\
\e(x,y)[\zeta(y),[z,x]\,]&=&\e(x,y)\mu\big(\zeta(y),\mu(z,x)\big)-\e(y,z)\mu\big(\mu(z,x),\zeta(y)\big)\\
&-&\e(x,y)\e(z,x)\mu\big(\zeta(y),\mu(x,z)\big)+\e(z,x)\e(y,z)\mu\big(\mu(x,z),\zeta(y)\big),\\
\e(y,z)[\zeta(z),[x,y]\,]&=&\e(y,z)\mu\big(\zeta(z),\mu(x,y)\big)-\e(z,x)\mu\big(\mu(x,y),\zeta(z)\big)\\
&-&\e(y,z)\e(x,y)\mu\big(\zeta(z),\mu(y,x)\big)+\e(x,y)\e(z,x)\mu\big(\mu(y,x),\zeta(z)\big).
\end{eqnarray*}
Then from Definition \ref{def} it follows that
$\circlearrowleft_{x,y,z}\e(z,x)[\zeta(x),[y,z]\,]=0.$
 \QED

By a homomorphism of Hom-Lie color algebras $ f :
\big(L,[\cdot,\cdot],\zeta,\e\big)\rightarrow
\big(L',[\cdot,\cdot]',\zeta',\e\big)$ we mean an algebra
homomorphism from $L$ to $L'$ such that the diagram below

\[\begin{CD}
L@>f>>L'\\
@V\zeta VV      @VV\zeta'V\\
L@>>f>L'\\
\end{CD}\]
commutes, i.e., $f\circ\zeta =\zeta'\circ f$.

 We now have the following result, which gives a way to construct
 Hom-Lie color algebras via a Lie color algebra together with an
 endomorphism.
\begin{theo}\label{th2} Let $\big(L, [\cdot,\cdot],\varepsilon\big) $ be a Lie color algebra and $\zeta: L\rightarrow
L$ be an even algebra endomorphism. Then
$\big(L,[\cdot,\cdot]_\zeta, \zeta,\varepsilon\big)$, where
$[x,y]_\zeta=\zeta\big([x,y]\big)$, is a Hom-Lie color algebra.
Moreover, suppose that $\big(L',[\cdot,\cdot]',\varepsilon\big)$ is
another Lie color algebra and $\zeta': L'\rightarrow L'$ be an even
algebra endomorphism. If $f:L\rightarrow L'$ is a Lie color algebra
homomorphism that satisfies $f\circ \zeta=\zeta'\circ f$, then
$$f:\big(L,[\cdot,\cdot]_{\zeta},\zeta,\varepsilon\big)\rightarrow
\big(L',[\cdot,\cdot]'_{ \zeta'},\zeta',\varepsilon\big)$$ is also a
homomorphism of Hom-Lie color algebras.
\end{theo}
\noindent{\it Proof}.\ \ It only needs to show that
$\big(L,[\cdot,\cdot]_\zeta, \zeta,\varepsilon\big)$ satisfies the
Hom-$\e$-Jacobi identity. Indeed, for all $x,y,z\in \mathcal
{H}(L)$, one has
\begin{eqnarray*}
\circlearrowleft_{x,y,z}\varepsilon(z,x)[\zeta(x),[y,z]_{\zeta}]_\zeta\,&=&\circlearrowleft_{x,y,z}\varepsilon(z,x)\zeta\Big([\zeta(x),\zeta\big([y,z]\big)]\Big)\\
&=&\circlearrowleft_{x,y,z}\varepsilon(z,x)\zeta^2\big([x,[y,z]\,]\big)\\
&=&\zeta^2\Big(\circlearrowleft_{x,y,z}\varepsilon(z,x)[x,[y,z]\,]\Big)\\
&=&0.
\end{eqnarray*}
Hence, $\big(L,[\cdot,\cdot]_\zeta, \zeta,\varepsilon\big)$ is a
Hom-Lie color algebra. The second assertion follows from
$$f\big([x,y]_\zeta\big)=f\circ \zeta\big([x,y]\big)=\zeta'\circ f\big([x,y]\big)=\zeta'\big([f(x),f(y)]'\big)=[f(x),f(y)]'_{\zeta'},$$
since $f$ is a Lie color algebra homomorphism satisfying $f\circ
\zeta=\zeta'\circ f$. \QED

\begin{ex}\label{example} \rm ({\bf Hom-Lie color $sl(2,\F)$}) \rm Assume that $sl(2,\F)$
is the three dimensional simple Lie algebra with the standard basis
\[
e=\left(\begin{array}{ccc} 0&1\\
0&0
\end{array}\right),\ \ \
f=\left(\begin{array}{ccc} 0&0\\
1&0
\end{array}\right),\ \ \
h=\left(\begin{array}{ccc} 1&0\\
0&-1
\end{array}\right),\ \ \
\]
satisfying
\begin{eqnarray*}\label{7}
[h,e]=2e,\ \ \ [h,f]=-2f,\ \ \ [e,f]=h.
\end{eqnarray*}
Put $a_1=\frac{i}{2}(e-f)$, $a_2=-\frac{1}{2}(e+f)$,
$a_3=\frac{i}{2}h$, where $i^2=-1$. Thus, we have
$$[a_1,a_2]=-a_3,\ \ [a_2,a_3]=a_1, \ \ [a_3,a_1]=a_2.$$
Let $\G=\Z_2\times\Z_2$. Then $sl(2,\F)=\oplus_{\a\in\G}\F X_{\a}$
is $\G$-graded with
$$X_{(0,0)}=0,\ \ X_{(1,0)}=a_1,\ \ X_{(0,1)}=a_2,\ \ X_{(1,1)}=a_3.$$
Define a mapping $\e:\G\times \G\rightarrow \F^*$ by
$$\e\big((\alpha_1,\alpha_2),(\beta_1,\beta_2)\big)=(-1)^{{\alpha_1\beta_1}+{\alpha_2\beta_2}},\ \ \forall \ \ \alpha_i,\beta_i\in\Z_2,\ i=1,2.$$
One can check that $\e$ is a  bi-character on $\G$. If the
$\G$-graded multiplication $\langle\cdot,\cdot\rangle$ turns
$sl(2,\F)$ into a Lie color algebra, then we have
\begin{eqnarray*}
\langle a_1,a_2\rangle=c_{12}a_3,\ \ \ \langle
a_2,a_3\rangle=c_{23}a_1,\ \ \ \langle a_3,a_1\rangle=c_{31}a_2.
\end{eqnarray*}
where $c_{12}$, $c_{23}$ and $c_{31}$ are arbitrary scalars in $\F$.
Now put $c_{12}=-1$, $c_{23}=1$, $c_{31}=1$. The Lie color algebra
$sl(2,\F)$ so defined is called the color analogue of the classical
Lie algebra $sl(2,\F)$.

Consider the nonzero linear map $\zeta: sl(2,\F)\rightarrow
sl(2,\F)$ defined by
\begin{eqnarray}\label{8}
\zeta(a_1)=x a_1,\ \ \ \zeta(a_2)=y a_2, \ \ \ \zeta(a_3)=z a_3,
\end{eqnarray}
where $x,y,z\in\F.$ It is not difficult to check that $\zeta$ is an
even endomorphism of Lie color algebra $sl(2,\F)$ if and only if
\begin{eqnarray*}
xy=z,\ \ \ yz=x, \ \ \ zx=y,
\end{eqnarray*}
which implies $x=y=z=\pm1.$ Now put $x=y=-1$ and $z=1$ in (\ref{8}).
Then according to Theorem \ref{th2}, we obtain a Hom-Lie color
algebra
$\big(sl(2,\F),\langle\cdot,\cdot\rangle_\zeta,\zeta,\e\big)$ with
the bracket $\langle\cdot,\cdot\rangle_\zeta$ on the basis
$\{a_1,a_2,a_3\}$ given by

\begin{eqnarray}\label{88}
 \langle a_1, a_2\rangle_\zeta=-a_3,\ \ \ \langle a_2, a_3\rangle_\zeta
=-a_1,\ \ \ \langle a_3, a_1\rangle_\zeta=-a_2.
\end{eqnarray}
\end{ex}

\begin{ex}\rm ({\bf Heisenberg Hom-Lie color algebra}) Let $H$ be the
three dimensional Heisenberg Lie algebra, which consists of the
strictly upper-triangular complex $3\times 3$ matrices. It has a
standard linear basis
\[
e_1=\left(\begin{array}{ccc} 0&1&0\\
0&0&0\\
0&0&0
\end{array}\right),\ \ \
e_2=\left(\begin{array}{ccc} 0&0&0\\
0&0&1\\
0&0&0
\end{array}\right),\ \ \
e_3=\left(\begin{array}{ccc} 0&0&1\\
0&0&0\\
0&0&0
\end{array}\right),\ \ \
\]
satisfying $$[e_1,e_2]=e_3,\  [e_1,e_3]=[e_2,e_3]=0.$$ Let
$\G=\Z_2\times\Z_2\times\Z_2$ with a bi-character $\e$ given by
\begin{eqnarray}\label{9}
\e(\alpha,\beta)=(-1)^{\alpha_1\beta_1+\alpha_2\beta_2+\alpha_3\beta_3},
\end{eqnarray}
for all $\a=(\a_1,\a_2,a_3)$ and $\b=(\b_1,\b_2,\b_3)$ in $\G$. Now
take $H$ to be a $\G $-graded linear space $H=\oplus_{\a\in\G}\F
X_{\a}$ with a homogeneous basis
$$ e_1=X_{(1,1,0)},\ \ e_2=X_{(1,0,1)},\
\ e_3=X_{(0,1,1)}.$$ The homogeneous components graded by the
elements of $\G$ different from $(1,1,0),\ (1,0,1)$ and $ (0,1,1)$
are zero and so are omitted. If the $\G$ graded multiplication
$\langle\cdot,\cdot\rangle$ turns $H$ into a Lie color algebra, then
with the bi-character $\e$ defined by (\ref{9}), we have
\begin{eqnarray*}\label{10}
\langle e_1,e_2\rangle=c_{12}e_3,\ \ \langle
e_2,e_3\rangle=c_{23}e_1,\ \ \langle e_3,e_1\rangle=c_{31}e_2,
\end{eqnarray*}
where $c_{12}$, $c_{23}$ and $c_{31}$ are arbitrary scalars in $\F$.
When $a$ and $b$ are in different homogeneous subspaces it follows
that $\langle a,b\rangle=\langle b,a\rangle$, whereas $\langle
a,b\rangle=-\langle b,a\rangle$ if $a$ and $b$ belongs to the same
one. Now put $c_{12}=1, \ c_{23}=0$ and $c_{31}=0$. The $\G$-graded
$\e$ Lie algebra $H$ so defined is the color analogue of the
Heisenberg Lie algebra, which is called {\it Heisenberg Lie color
algebra}.\\

Let $\lambda_1$ and  $\lambda_2$ be nonzero scalars in $\F$.
Consider the map $\zeta_{\lambda_1,\lambda_2}:H\rightarrow H$
defined on the basis elements by
\begin{eqnarray*}\label{10}
\zeta_{\lambda_1,\lambda_2}(e_1)=\lambda_1e_1,\ \
\zeta_{\lambda_1,\lambda_2}(e_2)=\lambda_2e_2,\ \
\zeta_{\lambda_1,\lambda_2}(e_3)=\lambda_1\lambda_2e_3.
\end{eqnarray*}
It is straightforward to check that $\zeta_{\lambda_1,\lambda_2}$
defines a Lie color algebra homomorphism. By virtue of Theorem
\ref{th2} above we have a Hom-Lie color algebra
$H_{\lambda_1,\lambda_2}=\big(H,\langle\cdot,\cdot\rangle_{\zeta_{\lambda_1,\lambda_2}},\zeta_{\lambda_1,\lambda_2},\e\big)$,
whose bracket satisfies the following non-vanishing twisted
Heisenberg relation
\begin{eqnarray*}
\langle
e_1,e_2\rangle_{\zeta_{\lambda_1,\lambda_2}}=\lambda_1\lambda_2e_3.
\end{eqnarray*}
We view the collection
$\{H_{\lambda_1,\lambda_2}:\lambda_1,\lambda_2\in\F^*\}$ as
two-parameter family of deformations of Heisenberg Lie color algebra
$H$ into Hom-Lie color algebras.
\end{ex}

 \begin{ex}\rm ({\bf Hom-Lie color algebra of Witt type}) Let
$s:\Gamma\rightarrow \F$ be a function defined over $\Gamma$ and
consider the $\Gamma$-graded vector space $L =\oplus_{\a\in\Gamma}
\mathbb{F}e_{\a}$ over $\F $ with a basis $\{e_{\a} |\a\in\Gamma\}$.
Assume that $\e:\Gamma\times\G\rightarrow \F^*$ is a bi-character on
$\Gamma$. The linear space $L$ endowed with the bracket defined by
$$[e_{\a},e_{\b}]= \big(s(\b)-\e(\a,\b)s(\a)\big)e_{\a+\b}$$
can be shown to be a Lie color algebra $ L=L(\Gamma,\e, s)$ under
certain conditions on the grading group $\Gamma$ and the mappings $
\e$ and $s$. Then $ L=L(\Gamma,\e, s)$ is said to be a Lie color
algebra of Witt type. For further details, we refer to [\ref{Z}].
Consider $L'=L\otimes \F[t,t^{-1}]$, where $\F[t,t^{-1}]$ is the
Laurent polynomials algebra. From Example \ref{ex1}, we know $L'$ is
also a Lie color algebra. Given any nonzero scalar $\lambda\in \F$,
we define a linear mapping $\zeta_\lambda: L'\rightarrow L'$ by
$\zeta_\lambda(e_{\a}\otimes f(t))=e_{\a}\otimes f(\lambda+t)$ for
any $e_\alpha\in L$, $f(t)\in \F[t,t^{-1}]$. Then $\zeta_\lambda$ is
a homomorphism of Lie color algebras. Indeed,
\begin{eqnarray*}
[\zeta_\lambda(e_{\a}\otimes f(t)),\zeta_\lambda(e_{\b}\otimes g(t))]_{L'} &=&[e_{\a}\otimes f(\lambda+t),e_{\b}\otimes g(\lambda+t)]_{L'}\\
&=&[e(\alpha),e(\beta)]_{L}\otimes f(\lambda+t)g(\lambda+t)\\
&=&\big(s(\b)-\e(\a,\b)s(\a)\big)e_{\a+\b}\otimes f(\lambda+t)g(\lambda+t) \\
&=&\big(s(\b)-\e(\a,\b)s(\a)\big)\zeta_\lambda\big(e_{\a+\b}\otimes f(t)g(t)\big)\\
&=&\zeta_\lambda\big([e_\a,e_\b]_L\otimes f(t)g(t)\big)\\
&=&\zeta_\lambda\big([e_{\a}\otimes f(t),e_{\b}\otimes
g(t)]_{L'}\big).
\end{eqnarray*}
Then
$L'_\lambda=\big(L',[\cdot,\cdot]_{\zeta_\lambda},\zeta_\lambda,\varepsilon\big)$
is a Hom-Lie color algebra thanks to Theorem \ref{th2} above. We
regard the collection $\{L'_\lambda: \lambda\in\F$\} as a
one-parameter family of deformations of the Lie color algebras of
Witt type into Hom-Lie color algebras.
\end{ex}
\vs{6pt}

\cl{\bf\S4. \ $\sigma$-twist of Hom-Lie color
algebra}\setcounter{section}{4}\setcounter{equation}{0}
 \vs{6pt}
 In this section, we shall establish a close relationship between
Hom-Lie color algebras corresponding to different form $\sigma$ on
$\Gamma.$ We provide conditions on $\sigma$ which ensure that the
$\sigma$-twist of Hom-Lie color algebra is also a Hom-Lie color
algebra. Finally, we give an example of $\sigma$-twist of Hom-Lie
color $sl(2,\F)$.
\par Let $\big(L,[\cdot,\cdot],\zeta,\e\big)$ be a Hom-Lie color
algebra. Given any mapping $\sigma:\Gamma\times\Gamma \rightarrow
\mathbb{F}^*$, we define on the $\Gamma$-graded vector space $L$ a
new multiplication $[\cdot,\cdot]^{\sigma}$ by the requirement that
\begin{eqnarray}\label{2}
[x,y]^{\sigma}=\sigma(x,y)[x,y],
\end{eqnarray}
 for all the homogeneous elements
$x,y$ in $L$. The $\Gamma$-graded vector space $L$, endowed with the
multiplication $[\cdot,\cdot]^{\sigma}$, is a $\Gamma$-graded
algebra which will be called a {\it $\sigma$-twist }of $L$ and will
be denoted by $L^{\sigma}$. We are now looking for conditions on
$\sigma$ which ensure that
$(L^{\sigma},[\cdot,\cdot]^{\sigma},\zeta,\e)$ is also a Hom-Lie
color algebra. \par It is easy to see that the bilinear mapping
$[\cdot,\cdot]^{\sigma}$ is $\e$ skew-symmetric if and only if\par
${\rm(I)}$ $\sigma$ is symmetric, i.e. $\sigma
(\b,\gamma)=\sigma(\gamma,\beta)$, for any
$\beta,\gamma\in\Gamma.$\\
Furthermore, the product $[\cdot,\cdot]^{\sigma}$ satisfies the
Hom-$\e$-Jacobi identity if and only if\par ${\rm(II)}$
$\sigma(\alpha,\beta)\sigma(\gamma,\alpha+\beta)$ is invariant under
cyclic permutations of $\alpha,\beta,\gamma\in\Gamma.$\par

We call such a mapping $\sigma:\Gamma\times\Gamma \rightarrow
\mathbb{F}^*$ satisfying both ${\rm(I)}$ and ${\rm(II)}$ is a {\it
symmetric multiplier} on $\G.$ Now from the discussion above it
follows:
\begin{prop}\label{p3} With notations above.
 Let $\big(L,[\cdot,\cdot],\zeta,\e\big)$ be a Hom-Lie color algebra and
 suppose that
 $\sigma$ is a symmetric multiplier on $\G$. Then the $\sigma$-twist $\big(L^{\sigma},[\cdot,\cdot]^{\sigma},\zeta,\e\big)$ is also a Hom-Lie
 color algebra under the same twisting map $\zeta$.
\end{prop}
\begin{coro}\label{co1}
Let $\big(L',[\cdot,\cdot]',\zeta',\e\big)$ be a second Hom-Lie
color algebra and $\sigma$ be a symmetric multiplier on $\G$. If
$f:L\rightarrow L'$ is a homomorphism of Hom-Lie color algebras,
then $f$ is also a homomorphism of Hom-Lie color algebras
$\big(L^{\sigma},[\cdot,\cdot]^{\sigma},\zeta,\e\big)$ into
$\big(L'^{\sigma},[\cdot,\cdot]'^{\sigma},\zeta',\e\big)$.
\end{coro}
\noindent{\it Proof}.\ \ One has $f\circ\zeta=\zeta'\circ f$ by the
assumption. For any homogenous element $x,y\in L$, we have
\begin{eqnarray*}
[f(x),f(y)]'^{\sigma}=\sigma\big(f(x),f(y)\big)[f(x),f(y)]'=\sigma(x,y)f\big([x,y]\big)=f\big(\sigma(x,y)[x,y]\big)=f\big([x,y]^{\sigma}\big),
\end{eqnarray*}
which proves the result.\QED\\
\begin{rema}\rm It is easy to construct a large class of
symmetric multipliers on $\G$ as follows. Let $\omega$ be an
arbitrary mapping of $\G$ into $\mathbb{F}^*$. Then the mapping
$\tau:\G\times\G\rightarrow \mathbb{F}^*$ defined by
$$\tau(\alpha,\beta)=\omega(\alpha+\beta)\omega(\alpha)^{-1}\omega(\beta)^{-1},\ \ \forall\  \alpha,\beta\in\G,$$
is a symmetric multiplier on $\G$.
\end{rema}

 Now given any mapping
$\sigma:\Gamma\times\Gamma \rightarrow \mathbb{F}^*$, which endows
$L$ with a new multiplication defined by (\ref{2}). we define a
mapping $\delta:\Gamma\times\Gamma \rightarrow \mathbb{F}^*$ by
\begin{eqnarray}\label{3}
\delta(\alpha,\beta)=\sigma(\alpha,\beta)\sigma(\beta,\alpha)^{-1},\
\ \mbox{for\ all}\ \alpha,\ \beta\in\Gamma.
\end{eqnarray}
Then it follows for any homogeneous elements $x,y\in L$ that
\begin{eqnarray}\label{4}
[x,y]^{\sigma}=-\e(x,y)\delta(x,y)[y,x]^{\sigma},
\end{eqnarray}
where we also simply write $\delta(x,y)$ instead of
$\delta(\bar{x},\bar{y})$.\par

 M. Scheunert in \cite{M1} provided the necessary and sufficient
 conditions on $\sigma$ for which ensure that $\varepsilon\delta$ is
 a bi-character on $\G$ (with $\e\delta(\alpha,\beta)=\varepsilon(\alpha,\beta)\delta(\alpha,\beta)$ for $\alpha,\beta\in\G$) and that $L^{\sigma}$ is an
 $\varepsilon\delta$-Lie algebra. Namely,
 $\big(L^{\sigma},[\cdot,\cdot]^{\sigma},\e\delta\big)$ is a $\G$-graded
 $\e\delta$-Lie algebra if and only if
 \begin{eqnarray}\label{5}
\sigma(\alpha,\beta+\gamma)\sigma(\beta,\gamma)=\sigma(\alpha,\beta)\sigma(\alpha+\beta,\gamma),\
\ \forall\ \ \a, \beta, \gamma\in \G.
\end{eqnarray}
\begin{defi} (See Definition 5 in \cite{M1}) Any mapping $\sigma:\Gamma\times\Gamma\rightarrow \mathbb{F}^*$
satisfying condition (\ref{5}) is called a {\it multiplier} on
$\Gamma.$
\end{defi}
\begin{rema}\rm Any symmetric multiplier that we have defined on $\G$ is a multiplier. That is why it is so named there.
\end{rema}

 For any multiplier $\sigma$, the mapping $\delta$ defined by
equation (\ref{3}) is a bi-character on $\G$ which is said to be
associated with $\sigma$. Note that $\delta(\a,\a)=1$ and it follows
from equation (\ref{5}) that $\sigma(0,\a)=\sigma(\a,0)=\sigma(0,0)$
for all $\a\in\G.$
\begin{prop}\label{p4} Let $\big(L,[\cdot,\cdot],\zeta,\e\big)$ be a Hom-Lie
color algebra. Suppose we are given a multiplier $\sigma$ on $\G$;
let $\delta$ be the bi-character on $\G$ associated with it. Then
$\big(L^{\sigma},[\cdot,\cdot]^{\sigma}, \zeta,
\varepsilon\delta\big)$ is a Hom-Lie color algebra.
\end{prop}
\noindent{\it Proof}.\ \ It only needs to verify the
Hom-$\e\delta$-Jacobi identity since the $\e\delta$-skew-symmetry
follows from equation (\ref{4}). For any homogeneous elements
$x,y,z\in L$, one has
\begin{eqnarray*}
\circlearrowleft_{x,y,z}\e\delta(z,x)[\zeta(x),[y,z]^{\sigma}]^{\sigma}=\sigma(x,y)\sigma(y,z)\sigma(z,x)\circlearrowleft_{x,y,z}\e(z,x)[\zeta(x),[y,z]\,]=0,
\end{eqnarray*}
 by the
assumption that $\big(L,[\cdot,\cdot],\zeta,\e\big)$ is a Hom-Lie
color algebra.\QED
\begin{coro}\label{co2}
Let $\big(L',[\cdot,\cdot]',\zeta',\e\big)$ be a second Hom-Lie
color algebra. Given a multiplier $\sigma$ on $\G$; let $\delta$ be
the bi-character on $\G$ associated with it. If $f:L\rightarrow L'$
is a homomorphism of Hom-Lie color algebras, then $f$ is also a
homomorphism of Hom-Lie color algebra
$\big(L^{\sigma},[\cdot,\cdot]^{\sigma},\zeta,\e\delta\big)$ into
Hom-Lie color algebra
$\big(L'^{\sigma},[\cdot,\cdot]'^{\sigma},\zeta,\e\delta\big)$.
\end{coro}\par

Now let us consider an example of $\sigma$-twist of Hom-Lie color
algebra.
\begin{ex} \rm We apply Proposition \ref{p4} to the Hom-Lie color algebra $\big(sl(2,\F), \langle\cdot,\cdot\rangle_\zeta, \zeta,\varepsilon\big)$
constructed in Example \ref{example}. Recall that
$\G=\Z_2\times\Z_2$ with the bi-character $\varepsilon$ defined by
$$\e\big((\alpha_1,\alpha_2),(\beta_1,\beta_2)\big)=(-1)^{{\alpha_{1}\beta_1}+{\alpha_2 \beta_2}},\ \ \forall \ \ \alpha_i,\beta_i\in\Z_2,\ i=1,2.$$
Define a mapping $\sigma: \G\times\G\rightarrow \F^*$ by
$$\sigma\big((\alpha_1,\alpha_2),(\beta_1,\beta_2)\big)=(-1)^{\alpha_1\beta_2},\ \ \forall \ \ \alpha_i,\beta_i\in\Z_2,\ i=1,2.$$
One can check that $\sigma$ satisfies condition (\ref{5}), namely,
it is a multiplier on $\G$. Then from (\ref{3}) we get a
bi-character $\delta$ by
$$\delta((\alpha_1,\alpha_2),(\beta_1,\beta_2))=(-1)^{{\alpha_1\beta_2}-{\alpha_2\beta_1}},\ \ \forall \ \ \alpha_i,\beta_i\in\Z_2,\ i=1,2.$$
Let $sl(2,\F)^{\sigma}$ be the $\sigma$-twist of the Hom-Lie color
algebra $\big(sl(2,\F), \langle\cdot,\cdot\rangle_\zeta,
\zeta,\varepsilon\big)$ with respect to $\sigma.$  It follows from
Proposition \ref{p4} that $\big(sl(2,\F)^{\sigma},
\langle\cdot,\cdot\rangle_\zeta^{\sigma},
\zeta,\varepsilon\delta\big)$ is also a Hom-Lie color algebra and
from (\ref{88}) and (\ref{2}) that the bracket
$\langle\cdot,\cdot\rangle_\zeta^\sigma$ can be explicitly given by
$$ \langle a_1, a_2\rangle_\zeta^\sigma=a_3,\ \ \ \langle a_2, a_3\rangle^\sigma_\zeta
=- a_1,\ \ \ \langle a_3, a_1\rangle_\zeta^\sigma=-a_2.
$$
\end{ex}

 \vs{6pt}

\cl{\bf\S5. \ Hom-Lie color admissible
algebra}\setcounter{section}{5}\setcounter{equation}{0}

 \vs{6pt}

 In this section, we aim to further extend the notions and results about Lie admissible
 algebras \cite{ME}
 to more generalized cases: Hom-Lie color admissible algebras and
flexible Hom-Lie color admissible algebras. We will also explore
some other general classes of such kind of algebras:
$G$-Hom-associative color algebras, using which we classify all the
Hom-Lie color admissible algebras.
\par  In this section, ``algebras'' means ``not necessarily associative
algebra''.
\begin{defi} Let $(A,\mu, \zeta)$ be a Hom-color algebra on the
$\G$-graded vector space $A$ defined by an even mutiplication $\mu$
and an algebra endomorphism $\zeta$. Let $\e$ be a bi-character on
$\Gamma$. Then $(A,\mu, \zeta)$ is said to be a hom-Lie color
admissible algebra if the bracket defined
 by
 \begin{eqnarray}
[x,y]=\mu(x,y)-\e(x,y)\mu(y,x)\label{6}
\end{eqnarray}
 satisfies the
Hom-$\e$-Jacobi identity
$$\circlearrowleft_{x,y,z}\e(z,x)[\zeta(x),[y,z]]=0,$$
for all the homogeneous elements $x,y,z\in A.$
\end{defi}
\begin{rema}\rm Since the color commutator bracket defined by (\ref{6}) is always
$\e$-skew symmetric, it makes any Hom-Lie color admissible algebra
into a
 Hom-Lie color algebra.\end{rema}

\begin{rema}\rm According to Proposition \ref{p1}, any associative Hom-color algebra $(L,\mu,\zeta)$ is a Hom-Lie  color admissible algebra with the same twisting
map $\zeta$. \end{rema}\par

Let $\big(L,[\cdot,\cdot],\zeta,\e\big)$ be a Hom-Lie color algebra.
Define a new commutator product $\langle\cdot,\cdot\rangle$ by
$$\langle x,y \rangle=[x,y]-\e(x,y)[y,x],\ \ \forall\ x,y\in\mathcal {H}(L).$$
 It is easy to see that
$\langle x,y\rangle=-\e(x,y)\langle y,x\rangle$. Moreover, we have
\begin{eqnarray*}
\circlearrowleft_{x,y,z}\e(z,x)\langle\zeta(x),\langle y,z\rangle\rangle&=&\circlearrowleft_{x,y,z}\e(z,x)\langle\zeta(x),[y,z]-\e(y,z)[z,y]\rangle\\&=&\circlearrowleft_{x,y,z}\e(z,x)\big(\langle\zeta(x),[y,z]\rangle-\e(y,z)\langle\zeta(x),[z,y]\rangle\big)\\
&=&\circlearrowleft_{x,y,z}\e(z,x)\Big([\zeta(x),[y,z]\,]-\e(x,y+z)[\,[y,z],\zeta(x)]\\&\
\ \
&-\e(y,z)[\zeta(x),[z,y]\,]+\e(y,z)\e(x,z+y)[\,[z,y],\zeta(x)]\Big)\\&=&
4\circlearrowleft_{x,y,z}\e(z,x)[\zeta(x),[y,z]\,]=0.
\end{eqnarray*}
Our discussion above now shows:
\begin{prop} Any Hom-Lie color algebra $\big(L,[\cdot,\cdot],\zeta,\e\big)$ is
 Hom-Lie color admissible.
\end{prop}\par
 Let $(L,\mu,\zeta)$
be a Hom-color algebra and $\e$ be a bi-character on the abelian
group $\G$. Let $$[x,y]=\mu(x,y)-\e(x,y)\mu(y,x), \ \mbox{for\ all}\
\  x,y\in \mathcal {H}(L),$$ be the associated color-commutator. An
{\it$\zeta$-associator} $a_{\mu,\zeta}$ of $\mu$ is defined by
\begin{eqnarray*}
a_{\mu,\zeta}(x,y,z)=\mu\big(\zeta(x),\mu(y,z)\big)-\mu\big(\mu(x,y),\zeta(z)\big),\
\ \forall\  x,y,z\in \mathcal {H}(L).
\end{eqnarray*}
 A Hom-color algebra
is said to be {\it flexible} if $a_{\mu,\zeta}(x,y,x)=0, \ \mbox{for
all} \ \ x,y\in \mathcal {H}(L).$\par  Now let us introduce the
notation:
\begin{eqnarray*}
S(x,y,z):=\e(z,x)a_{\mu,\zeta}(x,y,z)+\e(x,y)a_{\mu,\zeta}(y,z,x)+\e(y,z)a_{\mu,\zeta}(z,x,y).
\end{eqnarray*}
Then we have the following properties:
\begin{lemm}{\label{lem}}
$S(x,y,z)=\e(z,x)[\zeta(x),\mu(y,z)]+\e(x,y)[\zeta(y),\mu(z,x)]+\e(y,z)[\zeta(z),\mu(x,y)].$
\end{lemm}
\noindent{\it Proof}.\ \ The assertion follows expanding the
commutators on the right hand side:
\begin{eqnarray*}
&&\e(z,x)[\zeta(x),\mu(y,z)]+\e(x,y)[\zeta(y),\mu(z,x)]+\e(y,z)[\zeta(z),\mu(x,y)]\\
&&=\e(z,x)\mu\big(\zeta(x),\mu(y,z)\big)-\e(x,y)\mu\big(\mu(y,z),\zeta(x)\big)+\e(x,y)\mu\big(\zeta(y),\mu(z,x)\big)\\
&&-\e(y,z)\mu\big(\mu(z,x),\zeta(y)\big)+\e(y,z)\mu\big(\zeta(z),\mu(x,y)\big)-\e(z,x)\mu\big(\mu(x,y),\zeta(z)\big)\\
&&=\e(z,x)a_{\mu,\zeta}(x,y,z)+\e(x,y)a_{\mu,\zeta}(y,z,x)+\e(y,z)a_{\mu,\zeta}(z,x,y)\\
&&=S(x,y,z).
\end{eqnarray*}\QED
\begin{prop} A Hom-color algebra $(L,\mu,\zeta)$ is Hom-Lie color
admissible if and only if it satisfies
$$S(x,y,z)=\e(x,y)\e(y,z)\e(z,x)S(x,z,y),\ \ \forall \ x,y,z\in \mathcal {H}(L).$$
\end{prop}
\noindent{\it Proof}.\ \  From  Lemma \ref{lem} it follows
\begin{eqnarray*}
&&S(x,y,z)-\e(x,y)\e(y,z)\e(z,x)S(x,z,y)\\&&=\e(z,x)[\zeta(x),\mu(y,z)]+\e(x,y)[\zeta(y),\mu(z,x)]+\e(y,z)[\zeta(z),\mu(x,y)]\\
&&-\e(x,y)\e(y,z)\e(z,x)\Big(\e(y,x)[\zeta(x),\mu(z,y)]+\e(x,z)[\zeta(z),\mu(y,x)]+\e(z,y)[\zeta(y),\mu(x,z)]\Big)\\
&&=\circlearrowleft_{x,y,z}\e(z,x)[\zeta(x),\mu(y,z)-\e(y,z)\mu(z,y)]\\
&&=\circlearrowleft_{x,y,z}\e(z,x)[\zeta(x),[y,z]\,],
\end{eqnarray*}
which proves the result.\QED

In the following, we explore some other general classes of Hom-Lie
color admissible algebras, $G-$Hom-associative color algebras,
extending the class of Hom-associative algebras and Hom-associative
superalgebras. We will provide a classification of Hom-Lie color
admissible algebras using the symmetric group $S_3$, whereas it was
classified in \cite{AS} and \cite{FA} for the Hom-Lie and Hom-Lie
super cases, respectively.
\par

Let $S_3$ be the symmetric group generated by $\sigma_1=(1\ 2)$,
$\sigma_2=(2\ 3)$. Let $L=(L,\mu,\zeta)$ be a Hom-color algebra.
Suppose that $S_3$ acts on $L^{\times3}$ in the usual way, i.e.
$\sigma(x_1,x_2,x_3)=\big(x_{\sigma(1)},x_{\sigma(2)},x_{\sigma(3)}\big)$,
for any $\sigma\in S_3$ and $x_1,x_2,x_3\in \mathcal {H}(L).$ For
convenience, we introduce a notion of a parity of the transposition
$\sigma_i$ with $i\in\{1,2\}$ by setting
$$|\sigma_i(x_1,x_2,x_3)|=\e(x_i,x_{i+1}),\ \ \mbox{for\ all} \ x_1,x_2,x_3\in \mathcal {H}(L).$$
 It is natural to assume that the parity of the identity is $1$ and
for the composition $\sigma_i\circ \sigma_j$, it is defined by
\begin{eqnarray*}
|\sigma_i\circ
\sigma_j(x_1,x_2,x_3)|&=&|\sigma_j(x_1,x_2,x_3)|\cdot|\sigma_i\big(\sigma_j(x_1,x_2,x_3)\big)|\\
&=&|\sigma_j(x_1,x_2,x_3)|\cdot|\sigma_i\big(x_{\sigma_j(1)},x_{\sigma_j(2)},x_{\sigma_j(3)}\big)|.
\end{eqnarray*}
One can define by induction the parity for any composition. Hence we
have
\begin{eqnarray*}
|{\rm id}(x_1,x_2,x_3)|&=&1,\\
|\sigma_1(x_1,x_2,x_3)|&=&\e(x_1,x_2),\\
|\sigma_2(x_1,x_2,x_3)|&=&\e(x_2,x_3),\\
|\sigma_1\circ\sigma_2(x_1,x_2,x_3)|&=&\e(x_2,x_3)\e(x_1,x_3),\\
|\sigma_2\circ\sigma_1(x_1,x_2,x_3)|&=&\e(x_1,x_2)\e(x_1,x_3),\\
|\sigma_2\circ\sigma_1\circ\sigma_2(x_1,x_2,x_3)|&=&\e(x_2,x_3)\e(x_1,x_3)\e(x_1,x_2),
\end{eqnarray*}
 for any homogeneous elements $x_1,x_2,x_3$ in $L$. Let $sgn(\sigma)$ denote the signature of $\sigma\in
 S_3$. We have the following useful lemma:
\begin{lemm} \label{lem8}With notations above. A Hom-color algebra $(L,\mu,\zeta)$ is
Hom-Lie color admissible if the following condition holds
\begin{eqnarray*}
\mbox{$\sum_{\sigma\in
S_3}$}(-1)^{sgn(\sigma)}|\sigma(x_1,x_2,x_3)|a_{\mu,\zeta}\circ
\sigma(x_1,x_2,x_3)=0,\ \ \forall\ x_1,x_2,x_3\in\mathcal {H}(L).
\end{eqnarray*}
\end{lemm}
\noindent{\it Proof}.\ \  It only needs to verify the
Hom-$\e$-Jacobi identity. By straightforward calculation, the
associated color commutator satisfies
\begin{eqnarray*}
\circlearrowleft_{x_1,x_2,x_3}\e(x_3,x_1)[\zeta(x_1),[x_2,x_3]]=\e(x_3,x_1)\mbox{$\sum_{\sigma\in
S_3}$}(-1)^{sgn(\sigma)}|\sigma(x_1,x_2,x_3)|a_{\mu,\zeta}\circ
\sigma(x_1,x_2,x_3)=0.
\end{eqnarray*}
\QED \par Let $G$ be a subgroup of $S_3$, any Hom-color algebra
$(L,\mu,\zeta)$ is said to be $G$-{\it Hom-associative }if the
following equation holds:
\begin{eqnarray*}
\mbox{$\sum_{\sigma\in
G}$}(-1)^{sgn(\sigma)}|\sigma(x_1,x_2,x_3)|a_{\mu,\zeta}\circ
\sigma(x_1,x_2,x_3)=0, \ \ \forall\ x_1,x_2,x_3\in\mathcal {H}(L).
\end{eqnarray*}
 The following
result is a graded version of that obtained in \cite{AS}.
\begin{prop}
Let $G$ be a subgroup of the symmetric group $S_3$. Then any
$G$-Hom-associative color algebra $(L,\mu,\zeta)$ is Hom-Lie color
admissible.
\end{prop}
\noindent{\it Proof}.\ \ The $\e$-skew symmetry follows straightaway
from the definition. Assume that $G$ is a subgroup of $S_3$. Then
$S_3$ can be written as the disjoint union of the left cosets of
$G$. Say $S_3=\bigcup_{\sigma\in I}\sigma G$, with $I\subseteq S_3$
and for any $\sigma, \sigma'\in I$,
$$\sigma\neq \sigma'\Rightarrow \sigma G\cap \sigma'G=\emptyset.$$
Then one has
\begin{eqnarray*}
&&\mbox{$\sum_{\sigma\in
S_3}$}(-1)^{sgn(\sigma)}|\sigma(x_1,x_2,x_3)|a_{\mu,\zeta}\circ
\sigma(x_1,x_2,x_3)\\&&= \mbox{$\sum_{\tau\in
I}$}\mbox{$\sum_{\sigma\in \tau
G}$}(-1)^{sgn(\sigma)}|\sigma(x_1,x_2,x_3)|a_{\mu,\zeta}\circ
\sigma(x_1,x_2,x_3)=0,
\end{eqnarray*}
for any homogeneous elements $x_1,x_2,x_3$ in $L$. Hence the result
follows from Lemma \ref{lem8}.\QED \par Now we provide a
classification of the Hom-Lie color admissible algebras via
$G$-Hom-associative color algebras. The subgroups of $S_3$ are
$$G_1=\{{\rm id}\},\ \ G_2=\{{\rm id},\sigma_1\},\ \ G_3=\{{\rm id},\sigma_2\},\ \ G_4=\{{\rm id},\sigma_2\sigma_1\sigma_2=(1\ 3)\},\ \ G_5=A_3,\ \ G_6=S_3,$$
where $A_3$ is the alternating subgroup of $S_3$. Then from Lemma
\ref{lem8} we can obtain the following types of Hom-Lie color
admissible
 algebras :

$\bullet$ $G_1-$Hom-associative color algebras are the Hom-color
algebras defined in Definition \ref{def};

$\bullet$ $G_2-$Hom-associative color algebras satisfy the condition
\begin{eqnarray}\label{11}
\mu\big(\zeta(x),\mu(y,z)\big)-\e(x,y)\mu\big(\zeta(y),\mu(x,z)\big)=\mu\big(\mu(x,y),\zeta(z)\big)-\e(x,y)\mu\big(\mu(y,x),\zeta(z)\big);
\end{eqnarray}

$\bullet$ $G_3-$Hom-associative color algebras satisfy the condition
\begin{eqnarray}\label{12}
\mu\big(\zeta(x),\mu(y,z)\big)-\e(y,z)\mu\big(\zeta(x),\mu(z,y)\big)=\mu\big(\mu(x,y),\zeta(z)\big)-\e(y,z)\mu\big(\mu(x,z),\zeta(y)\big);
\end{eqnarray}

$\bullet$ $G_4-$Hom-associative color algebras satisfy the condition

$$\mu\big(\zeta(x),\mu(y,z)\big)-\mu\big(\mu(x,y),\zeta(z)\big)=\e(x,y)\e(y,z)\e(x,z)\Big(\mu\big(\zeta(z),\mu(y,x)\big)-\mu\big(\mu(z,y),\zeta(x)\big)\Big);$$

$\bullet$ $G_5-$Hom-associative color algebras satisfy the condition
\begin{eqnarray*}
\mu\big(\zeta(x),\mu(y,z)\big)-\e(x,y+z)\mu\big(\zeta(y),\mu(z,x)\big)-\e(x+y,z)\mu\big(\zeta(z),\mu(x,y)\big)\
\ \ \\=
\mu\big(\mu(x,y),\zeta(z)\big)-\e(x,y+z)\mu\big(\mu(y,z),\zeta(x)\big)-\e(x+y,z)\mu\big(\mu(z,x),\zeta(y)\big);
\end{eqnarray*}

$\bullet$ $G_6-$Hom-associative color algebras are the Hom-Lie color
admissible algebras by Lemma \ref{lem8},\\where $x,y,z$ are
homogeneous elements.

\begin{rema}\rm
When $\zeta$ is the identity map and $\e$ is trivial on $\G$, the
algebra defined by equation (\ref{11})
 is the classical Vinberg algebra or left symmetric algebra, while
 that defined in (\ref{12})
is Pre-Lie algebra or right symmetric algebra. Hence we will call
the algebras given by (\ref{11}) and (\ref{12}) the Hom-Lie color
Vinberg algebra and Hom-pre-Lie color algebra, respectively.
 \end{rema}

\vspace{4mm} \noindent\bf{\footnotesize Acknowledgements}\quad\rm
 {\footnotesize
This work was supported by the National Natural Science Foundation
of
China (Grant No. 10825101).}\\[4mm]

\end{document}